%Paper: alg-geom/9411010
%From: Miles Reid <miles@maths.warwick.ac.uk>
%Date: Wed, 16 Nov 1994 14:00:39 GMT
%Date (revised): Wed, 10 Jan 1996 10:13:15 +0100 (MET)

%%%%%%%%%%%%%%%%%%%%%%%%%%%%%%%%%%%%%%%%%%%%%
% Uses amstex, amsppt and optional epsf.tex %
%%%%%%%%%%%%%%%%%%%%%%%%%%%%%%%%%%%%%%%%%%%%%
% This paper comes with 2 postscript figures contained in a .uu file.
% To insert the figures you must be on a UNIX system; save the
% tex file as McKay.tex and edit out the mail header, then perform the
% following steps: (a) to unpack the uuencoded files containing the figures:
% Step 1. Save the .uu file as figures.uu, and edit out the mail header.
% Step 2. Do    csh figures.uu -- this run figures.uu as a program, which
% unpacks 3 files epsf, cube.ps and ram.ps.
% (b) to process the paper together with the figures:
% Step 3. Do   tex McKay.tex -- this produces McKay.dvi which is
% knowlegeable about the two postscript files.
% Step 4. Do   dvips McKay.dvi  to translate the dvi into a postscript file
% McKay.ps, which you can send to a postscript printer by lpr as usual.
%%%%%%%%%%%%%%%%%%%%%%%%%%%%%%%%%%%%%%%%%%%%%
\input amstex
\documentstyle{amsppt}
\magnification\magstep1
\TagsOnRight
\input epsf \epsfverbosetrue

                 % union over set
                 % intersection over set
\define\1{^{-1}}                       % inverse
\define\third{\frac13}                 % \third(x,y,z) gives 1/3(x,y,z)
\define\rth{\frac1r}                   % \rth(x,y,z) gives 1/r(x,y,z)
\define\into{\hookrightarrow}          % injective arrow
       % surjective arrow
\define\Span#1{\left<#1\right>}        % <x> span or hull of x
\define\QED{\ifhmode\unskip\nobreak\fi\quad Q.E.D.}  % QED
\define\Box{\square}
\define\GL{\mathop{{\roman{GL}}}}      % full linear group
\define\BD{\mathop{{\roman{BD}}}}      % binary dihedral group BD_4n
\define\BT{\mathop{{\roman{BT}}}}      % binary tetrahedral group BT_48
\define\SL{\mathop{{\roman{SL}}}}      % special linear group
\define\SO{\mathop{{\roman{SO}}}}      % orthogonal group
\define\Ext{\mathop{{\roman{Ext}}}\nolimits} % Ext group
\define\Fix{\operatorname{Fix}}        % Fixed locus
\define\Gal{\mathop{{\roman{Gal}}}}    % Galois group
\define\Hom{\operatorname{Hom}}        % Hom group
\define\Ram{\operatorname{Ram}}        % Ramification group
\define\Gr{\operatorname{Gr}}          % associated graded
\define\Stab{\operatorname{Stab}}      % Stabiliser group
\define\Proj{\operatorname{Proj}}      % homogeneous spectrum of ring
\define\Spec{\mathop{{\roman{Spec}}}}  % spectrum of ring
 % Trace map
\define\diag{\operatorname{diag}}      % diagonal group
\define\id{\roman{id}}                 % identity map
\define\dd{\roman{d}}                  % roman differential dx, dy
\define\age{\operatorname{age}}
\define\weight{\operatorname{weight}}
\define\rest#1{_{\textstyle{\vert}#1}} % restriction of map to subset

% Bourbaki open-face letters
\define\C{\Bbb C} % complex numbers
 % quaternions
 % natural numbers
\define\Q{\Bbb Q} % rational numbers
\define\R{\Bbb R} % real numbers
\define\Z{\Bbb Z} % integers
 % scroll F or Galois field F_q
 % additive or mult groups G_a, G_m
 % affine space
\define\proj{\Bbb P} % projective space

% script
\define\Oh{{\Cal O}} % structure sheaf
\define\sF{{\Cal F}} % sheaf F

% bold Greek
\loadbold
\define\bmu{\boldsymbol\mu}

% short Greeks
\define\al{\alpha}
\define\be{\beta}

\define\ep{\varepsilon}
\define\fie{\varphi}

\define\la{\lambda}
\define\si{\sigma}

\define\Ga{\Gamma}
\define\De{\Delta}

\define\Si{\Sigma}

\define\bij{\mathbin{\;\longleftrightarrow\;}}
\define\iso{\cong}
\define\wave{\widetilde}

 \topmatter
 \author Yukari Ito and Miles Reid \endauthor
 \title The McKay correspondence for finite subgroups of $\SL(3,\C)$ \endtitle
 \rightheadtext{McKay for $\SL(3,\C)$}
 \abstract Let $G\subset\SL(n,\C)$ be a finite subgroup and $\fie\colon Y\to
X=\C^n/G$ any resolution of singularities of the quotient space. We prove that
crepant exceptional prime divisors of $Y$ correspond one-to-one with ``junior''
conjugacy classes of $G$. When $n=2$ this is a version of the McKay
correspondence (with irreducible representations of $G$ replaced by conjugacy
classes). In the case $n=3$, a resolution with $K_Y=0$ is known to exist by
work of Roan and others; we prove the existence of a basis of $H^*(Y,\Q)$ by
algebraic cycles in one-to-one correspondence with conjugacy classes of $G$.
Our treatment leaves lots of open problems.
 \endabstract
 \address Dept.\ of Math Sciences, Univ.\ of Tokyo, Tokyo 113, Japan
\endaddress
 \email Yukari-I\@Tansei.CC.U-Tokyo.Ac.Jp \endemail
 \address Math Inst., Univ.\ of Warwick, Coventry CV4 7AL, England \endaddress
 \email  Miles\@Maths.Warwick.Ac.UK \endemail
 \endtopmatter
 \document
\head Contents \endhead

{{\narrower\narrower

1. Statement of the results

2. Proofs

3. Examples

4. Discussion

References

}}

\head 1. Statement of the results \endhead

Let $G\subset\SL(n,\C)$ be a finite subgroup and $X=\C^n/G$ the quotient space,
an affine variety with $K_X=0$. A {\it crepant resolution} $f\colon Y\to X$ is
a
resolution of singularities such that $K_Y=f^*K_X=0$. A crepant resolution does
not necessarily exist in dimension $\ge4$ (see 4.5); it is known
to exist in dimension 2 (classical, Du Val), and in dimension 3 by work of a
number of people (see for example \cite{Markushevich}, \cite{Ito1--3} and
\cite{Roan}), but the proofs are computational rather than conceptual.

This paper contributes some very easy remarks to the following question raised
by \cite{Dixon--Harvey--Vafa--Witten}, worked out by
\cite{Hirzebruch--H\"ofer},
and now famous among algebraic geometers as the ``Physicists' Euler number
conjecture''. See for example \cite{Roan} for the background.

 \proclaim{Conjecture~1.1} $G\subset\SL(n,\C)$ is a finite subgroup, $X=\C^n/G$
the quotient space and $f\colon Y\to X$ a crepant resolution. Then there exists
a basis of\/ $H^*(Y,\Q)$ consisting of algebraic cycles in one-to-one
correspondence with conjugacy classes of\/ $G$.
 \endproclaim

It is an elementary fact that $Y$ has no odd-dimensional cohomology, and that
$H^{2i}(Y,\Q)$ is spanned by algebraic cycles (see 4.1).

 \definition{Definition~1.2} Let $G$ be a finite group. Write
$\Ga=\Hom(\bmu_R,G)$, where $R$ is any common multiple of the orders of all
$g\in G$ and $\bmu_R$ is the group of complex $R$th roots of $1$; thus
$\Ga$ is a set which becomes a group isomorphic to $G$ after making a choice
of roots of 1. In colourful language, you can think of $\Ga$ as the ``Tate
twist'' $\Ga=G(-1)=\Hom(\bmu,G)$, where $\bmu=\widehat\Z(1)=\varprojlim\bmu_r$.
 \enddefinition

 \proclaim{Theorem~1.3 \rm (age grading)} Let\/ $G\subset\SL(n,\C)$ be a finite
subgroup and\/ $\Ga$ as above. Then $\Ga$ has a canonical grading
 $$
\Ga=\sqcup_{i=1}^{n-1}\Ga_i,\qquad\text{with}\quad\age\Ga_i=i.
 $$
(Here $\Ga_0=\{e\}$, the identity element of\/ $\Ga$.)
 \endproclaim

Putting a grading on $G$ is a well-known idea in toric geometry; see for
example \cite{C3-f}, Theorem~3.1, and 2.2 below, and compare
\cite{Batyrev--Dais}. The point of the definition of $\Ga$ is to make the
grading independent of the choice of roots of 1. If you believe
Conjecture~1.1, you expect conjugacy classes in $\Ga_i$ to give the basis
of $H^{2i}$, so that the age grading of $\Ga$ is an analog of the usual
weight filtration of $H^*(Y)$, with $\weight=2\times\age$.

The {\it junior elements} of $\Ga$ (or of $G$) are the elements $g\in\Ga_1$ of
minimal age 1; the set $\Ga_1$ is called the {\it junior simplex} of $\Ga$ (or
of $G$). The grading of Theorem~1.3 is invariant under conjugacy by $G$ (or by
$\Ga$, it's the same thing). The notation $\Ga_i/G$ always means conjugacy
classes of elements of age $i$.

 \proclaim{Theorem~1.4} There is a canonical one-to-one correspondence between
junior conjugacy classes in $\Ga$ and crepant discrete valuations of $X$.
\endproclaim

The statement is explained in 2.3 below. Roughly, a discrete valuation of $X$
is a prime divisor $E\subset V$ on any normal model $V$ of $X$ up to birational
equivalence. Discrete valuations enable us to state Theorem~1.4 in a
birationally invariant language, without any assumptions on crepant resolutions
or minimal models of $X$.

A {\it minimal model\/} of $X$ (in the sense of Mori theory) is a projective
morphism $f\colon Y\to X$, where $f$ is crepant (that is, $K_Y=f^*K_X$) and $Y$
has $\Q$-factorial terminal singularities.

 \proclaim{Corollary~1.5} Suppose that there exists a minimal model of
$X=\C^n/G$. Then there are canonical one-to-one correspondences
 $$
\Ga_1/G
\bij
\text{exceptional prime divisors of $f$}
\bij
\text{basis of $H^2(Y,\Q)$.}
 $$
 \endproclaim

Now suppose that $n=3$. Then it is known by work of Mori and others (these
are general arguments, applicable to any 3-fold) that a minimal model
$f\colon Y\to X=\C^3/G$ exists. Moreover, it has factorial cDV singularities
(see \cite{Koll\'ar}, 2.1.7.6), and is therefore a $\Q$@-cohomology manifold,
so that Poincar\'e duality applies between $H^4(Y,\Q)$ and cohomology with
compact support $H^2_c(Y,\Q)$.

 \proclaim{Theorem~1.6} Let $f\colon Y\to X$ be as just described. Then
$H^2_c(Y,\Q)$ is based by the compactly supported exceptional prime divisors
of\/ $Y$, that is, the divisors $E\subset f\1(0)$ lying over $0$. Write
$\Ga_1^{(0)}$ for the junior elements of\/ $\Ga$ corresponding to these (see
also 2.10 below).

Then there are canonical one-to-one correspondences
 $$
\Ga_2/G
\bij
\Ga_1^{(0)}/G
\bij
\text{basis of $H^2_c(Y,\Q)$}
\bij
\text{dual basis of $H^4(Y,\Q)$}.
 $$
 \endproclaim

These results are proved in \S2. Some features and implications are discussed
in \S4.

 \subhead 1.7. Acknowledgements \endsubhead We are grateful to the organisers
of
the Trento conference (Jun 1994), where a seminar was held which discussed lots
of this stuff (see for example \cite{Batyrev}). A.~Corti pointed out that
3-folds with $\Q$-factorial terminal singularities are $\Q$-homology manifolds,
so that Poincar\'e duality holds. Questions by John Moody have helped us to
improve the argument in a number of places. The second author wishes to thank
S.~Mori and Kyoto Univ.\ Koenkai for generous invitation and financial support.
We are also indebted to David Morrison for pointing out that physicists are
human beings and have individual names.

\head 2. Proofs \endhead

 \subhead 2.1. Proof of Theorem~1.3 \endsubhead Let $G\subset\SL(n,\C)$ be a
finite group. Any element $g\in G$ has $n$ eigenvalues $\la_1,\dots,\la_n$,
which are $r$th roots of 1 if $g^r=1$. In the normal run of things, in order to
write $\la_i=\ep^{a_i}$, we have to choose $\ep$ a primitive $r$th root of 1.
However, the present set-up side-steps this problem by considering elements
$g\in\Ga=\Hom(\bmu_R,G)$; then $\ep\to g(\ep)\to\la_i(g(\ep))$ is an unordered
set of $n$ characters of $\bmu_R$, and thus $\la_i(g(\ep))=\ep^{a_i}$. If $g$
has order $r$ then $a_1,\dots,a_n$ correspond to the characters
$\ep\mapsto\ep^{a_i}$ of order $r$, and we write them as the fractional
expression $\rth(a_1,\dots,a_n)$ with integers
$a_i\in[0,r)$.

Now, since $G\subset\SL(n)$ it follows that $\sum a_i\equiv0\mod r$, and hence
$\rth\sum a_i=\al$ for some integer $\al=\age g\in[0,n)$. This defines
the age grading of $\Ga$, and so completes the proof. \QED

\remark{Remarks} 1. In fixed coordinates $x_1,\dots,x_n$ on $\C^n$, the
expression $g=\rth(a_1,\dots,a_n)$ is a well-defined homomorphism
$\bmu_r\to\GL(n,\C)$, namely
 $$
\ep\mapsto\diag(\ep^{a_1},\dots,\ep^{a_n}).
 $$

2. Note that we didn't need to make an explicit choice of eigencoordinates,
since the characters $\ep^{a_i}$ are obtained as the set of eigenvalues of
matrixes, and so are invariant under conjugacy by $\GL(n,\C)$. We chose an
order $\rth(a_1,\dots,a_n)$ of the characters for convenience, but the
order disappears at the end, since the age is the elementary symmetric
function $\rth\sum a_i$.

3. The group $G$ itself does not have a well-defined grading: for example, if
$\ep$ is a primitive 7th root of 1 and $G=\Z/7$ is generated by
$\diag(\ep,\ep^2,\ep^4)\in\SL(3,\C)$, then choosing the same $\ep$ gives
 $$
\age(\overline1,\overline2,\overline4)=1\quad\text{and}\quad
\age(\overline3,\overline5,\overline6)=2,
 $$
whereas the choice of $\eta=\ep\1$ gives the opposite grading.
\endremark

 \subhead 2.2. Notation: the unit box $\Box$ and junior simplex $\De_1$
\endsubhead We now go through the same proof in more conventional terms,
and, at the same time, introduce notation that we need later.

Let $G\subset\SL(n,\C)$ be as usual and $A\subset G$ an Abelian group;
write $r$ for some exponent of $A$. Choose coordinates $x_1,\dots,x_n$
on $\C^n$ to diagonalise the action, so that $A$ acts on the
$x_1,\dots,x_n$ by diagonal matrixes:
 $$
g\mapsto\diag(\la_1(g),\dots,\la_n(g)),
 $$
where $\la_i\colon A\to\bmu_r$ are characters of $A$. Then the Pontryagin dual
of $A\into\bmu_r^n\subset(\C^\times)^n$ is a surjective homomorphism
 $$
h\colon \Z^n\to A^\vee\qquad
\text{defined by }
\prod x_i^{m_i}\mapsto \sum m_i\la_i.
 $$
Here $A^\vee=\Hom(A,\C^\times)$ is the character group of $A$, and we
identify $\Z^n$ with the
set of Laurent monomials $x^m=\prod x_i^{m_i}$, as usual in toric geometry. The
kernel of $h$ is the set of invariant monomials $M\subset\Z^n$, a sublattice of
finite index. The dual lattice of $M$ is an overlattice $L\supset\Z^n$. For
example, if $A$ is cyclic with generator $\diag(\ep^{a_1},\dots,\ep^{a_n})$
then
$L=\Z^n+\Z\cdot\rth(a_1,\dots,a_n)$. To distinguish the two, $M$ is
sometimes called the {\it lattice of monomials}, and $L$ the {\it lattice of
weights}.

Now in the copy of $\R^n$ corresponding to $\Z^n\subset L$, write
 $$
\Box=\prod[0,1)=\bigl\{(\al_1,\dots,\al_n)\in\R^n\,\bigm|\,0\le\al_i<1\bigr\}
 $$
for the unit cube, and
 $$
\De_1=\bigl\{(\al_1,\dots,\al_n)\in\Box\,\bigm|\,\sum\al_i=1\bigr\}
 $$
for the unit simplex (see Figure~1). The point is that, provided roots of 1 are
chosen, the quotient $L/\Z^n$ is naturally isomorphic to $A$, and since $\Box$
is a fundamental domain for the action of $\Z^n$ on $\R^n$, every element of
$L/\Z^n\iso A$ has a unique representative in $\Box$, and so $A=L\cap\Box$.
Thus the age grading of $A$ is just the natural slicing of $\Box$ by the planes
$\sum\al_i=j$ for $j=1,\dots,n-1$. In particular the set of junior elements of
$A$ equals $L\cap\De_1$.

\midinsert
\noindent
\hfil
\epsfbox{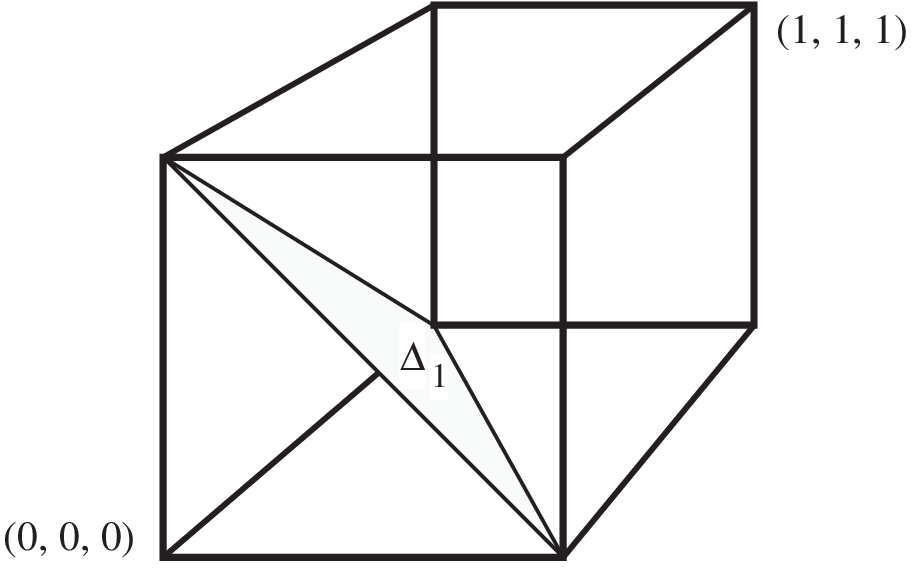}
\hfil
\par\vskip.2cm
\noindent
\hfil
{\bf Figure~1.} The unit box $\Box$
and junior simplex $\De_1$
\hfil
\endinsert

 \proclaim{Lemma} There is a canonical isomorphism
$L/\Z^n=A(1)=\Hom(\bmu_r,A)$.
 \endproclaim

 \demo{Proof} Given the coordinate system $x_1,\dots,x_n$, the exact sequence
 $$
0\to M\to\Z^n\to A^\vee\to0
 $$
is intrinsic, and, at first sight, it seems that the isomorphism
$L/\Z^n\to A$ should also be: in fact, lattice theory gives a canonical
duality pairing $\Z^n/M\times L/\Z^n\to\Q/\Z$. The subtle point is that
this is duality with coefficients in $\Q/\Z$, giving
 $$
L/\Z^n=\Hom(A^\vee,\Q/Z),
 $$
whereas the character group $A^\vee=\Hom(A,\C^\times)$ is the dual
with coefficients in $\C^\times$. Now obviously
 $$
A(1)=\Hom(\bmu_r,A)=\Hom(A^\vee,\Q/\Z).
 $$

An alternative, slightly bizarre, proof is to apply
$\Hom_\Z(\text{blank},\Z)$ to the exact sequence $0\to M\to\Z^n\to
A^\vee\to0$, to get
 $$
0=\Hom(A^\vee,\Z)\to\Z^n\to L\to\Ext^1_\Z(A^\vee,\Z)\to0,
\quad\text{that is,}\quad L/\Z^n=\Ext^1_\Z(A^\vee,\Z),
 $$
and finally to argue that
 $$
\Ext^1(A^\vee,\Z)=A\otimes\Ext^1(\widehat\Z,\Z)=A\otimes\widehat\Z(1).
\roman{\QED}
 $$
\enddemo

 \subhead 2.3. Overview of proof of Theorem~1.4 \endsubhead We start by
explaining the statement and the strategy. Let $f\colon V\to X$ be any
resolution of singularities, with exceptional divisors $\{E\}$. Write
$K_V=f^*K_X+\sum a_EE$. The discrepancy $a_E$ is independent of $V$. It is
known that $a_E\ge0$ for all $E$ (see 2.4 below, and compare \cite{C3-f},
Theorem~3.1). The {\it crepant} divisors $E$ are those with $a_E=0$. Every
crepant $E$ must appear on any resolution $V$. The theorem states that the
crepant divisors $E$ for any $V$ correspond one-to-one with $\Ga_1/G$.

Let $V$ be a normal variety (in applications, a blowup of $\C^n$ or a partial
resolution of $X$). Counting the zeros and poles of rational functions $h\in
k(V)$ along a prime divisor $E\subset V$ defines a discrete valuation of the
function field $k(V)$, which is a homomorphism $v_E\colon k(V)^\times\to\Z$
satisfying a suitable compatibility with the additive structure of $k(V)$. The
theory of discrete valuations is a device for discussing properties of
$E\subset V$ in a birational way. In what follows, although not absolutely
essential, it is convenient to use the definition of discrete valuation $v$ of
$k(V)$, and the fact that every geometric discrete valuation of a function
field $k(X)$ is of the form $v_E$ for some prime divisor $E$ on a partial
resolution $V\to X$. ({\it Geometric} means that the residue field $k_v$ has
transcendence degree $n-1$; this is a very easy case of Zariski's ``local
uniformisation'', see Zariski and Samuel \cite{Z--S}, Vol.~II, Chap.\ VI,
\S14, Theorem~31.)

It turns out that the only exceptional divisors of partial resolutions $V\to X$
we need arise from weighted blowups of $\C^n$: for each $g\in G$, take
eigencoordinates $x_1,\dots,x_n$ for $g$, and a weighting $\be=(b_1,\dots,b_n)$
on $x_1,\dots,x_n$ (usually closely related to the eigenvalues of $g$; see 2.7
for more precise conditions). The exceptional divisor of the weighted blowup
$W_\be\to\C^n$ is a weighted projective space $F_\be=\proj(b_1,\dots,b_n)$, and
the valuation $v_\be=v_{F_\be}$ of $k(\C^n)$ along $F_\be\subset W_\be$ takes
$x_i$ to $b_i$. The corresponding exceptional divisor $E$ of $X=\C^n/G$ is
somewhat messy to describe geometrically (see 2.6 below), but the valuation
$v_E$ is just the restriction of $v_\be$ to the subfield $k(X)=k(\C^n)^G\subset
k(\C^n)$, divided by the ramification degree. A valuation of this form is
called
a {\it monomial valuation} of $k(X)$.

Theorem~1.4 is a consequence of two easy tricks: (I) if $X=\C^n/G$ is a
quotient singularity and $E\subset V\to X$ a prime divisor on a partial
resolution, {\it ramification theory} reduces the calculation of the
discrepancy of $E$ to the case of a cyclic subgroup $\Ram F=\Z/r\subset G$ (see
2.5--6). (II) A cyclic quotient singularity $\C^n/(\Z/r)$ is a toric variety,
so
has a toric resolution (see 2.4); every valuation appearing on it is monomial
in the sense just discussed. Therefore every nonmonomial valuation of $X$ is
discrepant.

The rest of the proof just boils down to calculating monomial valuations and
the effect of the group action on them.

 \subhead 2.4. The cyclic case \endsubhead Consider first the case of a cyclic
group $\Span g$ with $g\in\Ga$ an element of order $r$. Choose eigencoordinates
to make $g=\rth(a_1,\dots,a_n)$ as described in 2.2. The following toric
construction is standard. (For details, compare, for example, the proof of
\cite{C3-f}, Theorem~3.1 or \cite{YPG}, \S4.)

The barycentric subdivision of the first octant of $\R^n$ at
$\rth(a_1,\dots,a_n)$ (with respect to the two different lattices
$\Z^n$ and $L$) defines two blowups
 $$
\spreadmatrixlines{3pt}
\matrix
F&\subset&W&\to&\C^n\\
&&\downarrow&&\downarrow\\
E&\subset&V&\to&\C^n/\Span g.
\endmatrix
 $$
Here $W$ is the weighted blowup of $\C^n$ with weights $(a_1,\dots,a_n)$. This
has a single exceptional component $F$, which is a weighted projective space
$F=\proj(a_1,\dots,a_n)$, and $V=W/\Span g$. It is easy to see that $F$ has
discrepancy $a_F=\sum a_i-1=\al r-1$ with $\al=\age g$ (see \cite{YPG},
Prop.~4.8), and that $\Span g$ acts biregularly on $W$, fixing $F$ pointwise,
and multiplying a normal coor\-dinate to $F$ by a primitive $r$th root of 1.
Thus the cyclic cover $W\to V$ has ramification of order $r$ along $E$, and
therefore $E$ has discrepancy $\rth(a_F-(r-1))=\al-1$. Therefore
$a_E\ge0$, and $E$ is crepant if and only if $\al=\age g=1$, that is, $g$ is
junior.

In this case, the quotient $X$ is a toric variety, so has a toric resolution
$V_\Si\to X$, obtained by subdividing the first octant of $\R^n$ into a fan
$\Si$ of basic simplexes for the lattice
$L=\Z^n+\Z\cdot\rth(a_1,\dots,a_n)$. The exceptional divisors of $V_\Si$
correspond to the 1@-skeleton of $\Si$; if $e=\frac1{r'}(b_1,\dots,b_n)\in L$
(for some
$r'\mid r$) is in the positive octant then its discrepancy is $\frac1{r'}\sum
b_i-1$. Thus (again because $g\in\SL(n,\C)$, so that $\sum b_i\equiv0\mod r'$),
every exceptional divisor has discrepancy $\ge0$, and the crepant divisors
correspond to $\De_1\cap L$, that is, to the junior elements $\Ga_1$.

 \subhead 2.5. Ramification groups \endsubhead The ramification theory of
discrete valuations in a Galois field extension $K\subset L$ is very similar to
the factorisation of prime ideals in a Galois extension of number fields. (See
2.6 below for the birational geometric picture.) Let $K\subset L$ with
$G=\Gal(L/K)$ be a finite Galois extension of function fields. A discrete
valuation $v\colon L^\times\to\Z$ of $L$ restricts to a map $v\rest K\colon
K^\times\to Z$ which obviously satisfies all the conditions for a discrete
valuation, except that it only maps onto a subgroup $r\Z\subset\Z$ of finite
index. We call $r$ the {\it ramification degree} of $v$ and write
$w=\rth(v\rest K)$. Conversely, a discrete valuation $w$ of $K$ extends in
finitely many ways to discrete valuations $v$ of $L$. These extensions are
permuted transitively by the Galois group, which acts on the valuations of $L$
by
 $$
g\colon v\mapsto v\circ g\qquad\text{for $g\in\Gal(L/K)$.}
 $$

Corresponding to a valuation $v$ of $L$, there is a tower of subgroups
 $$
\Ram v\subset\Stab v\subset\Gal(L/K).
 $$
The {\it stabiliser group} $\Stab v$ and the {\it ramification group} $\Ram
v$ are defined by
 $$
 \Stab v=\bigl\{g\in G\,\bigm|\,v\circ g=v\bigr\}\quad\text{and}\quad
 \Ram v=\bigl\{g\in\Stab v\,\bigm|\,g^*=\id\text{ on }k(v)\bigr\},
 $$
where $k(v)=A_v/m_v$ is the residue field of $v$. (Compare \cite{Z--S},
Vol.~II, Chap.~VI, \S12; these are also traditionally called the {\it
decomposition subgroup} $G_Z$ (Zerlegungsgruppe) and the {\it inertia
subgroup} $G_T$ (Tr\"agheitsgruppe) respectively.)

A crucial point for us is that $\Ram v$ is a cyclic group. The complex analysis
argument is transparent and is given at the end of 2.6. For convenience, we
spell out the algebraic proof of Zariski and Samuel: in their notation, by
\cite{Z--S}, Vol.~II, Chap.~VI, \S12, Theorem~24, Corollary, we have $G_V=1$;
moreover (same book, pp.~75--76), $G_T/G_V$ is dual to the extension of the
value groups $\Ga^*/\Ga$, which for discrete valuations is a finite cyclic
group of order $r$, where $r$ is the ramification degree. In either case,
characteristic zero is used in an essential way.

 \subhead 2.6. The geometric picture \endsubhead For a finite subgroup
$G\subset\GL(n,\C)$, the set-up is as follows. Write $L=k(\C^n)$ and
$K=k(X)=L^G$. A discrete valuation $v$ of $L$ is of the form $v=v_F$, where
$F\subset W$ is a prime divisor on a model $W\to\C^n$. There is a partial
resolution $V\to X$ of $X$ such that $F$ maps generically onto a prime divisor
$E\subset V$. You can prove this by applying the local uniformisation result
quoted in 2.3. The alternative is to construct $V$ as a quotient of some blowup
of $W$ on which $G$ acts biregularly. For this, let $\fie_g\colon W^g\to\C^n$
be
the morphism with $W^g=W$ and $\fie_g$ the composite $W\to\C^n@>g>>\C^n$, and
define $W'$ to be the birational fibre product of all the $W^g$. In other
words, take the closure in the fibre product of the graph of the locus where
all the $g$ are regular, normalised if necessary.

\midinsert
\noindent
\hfil
\epsfbox{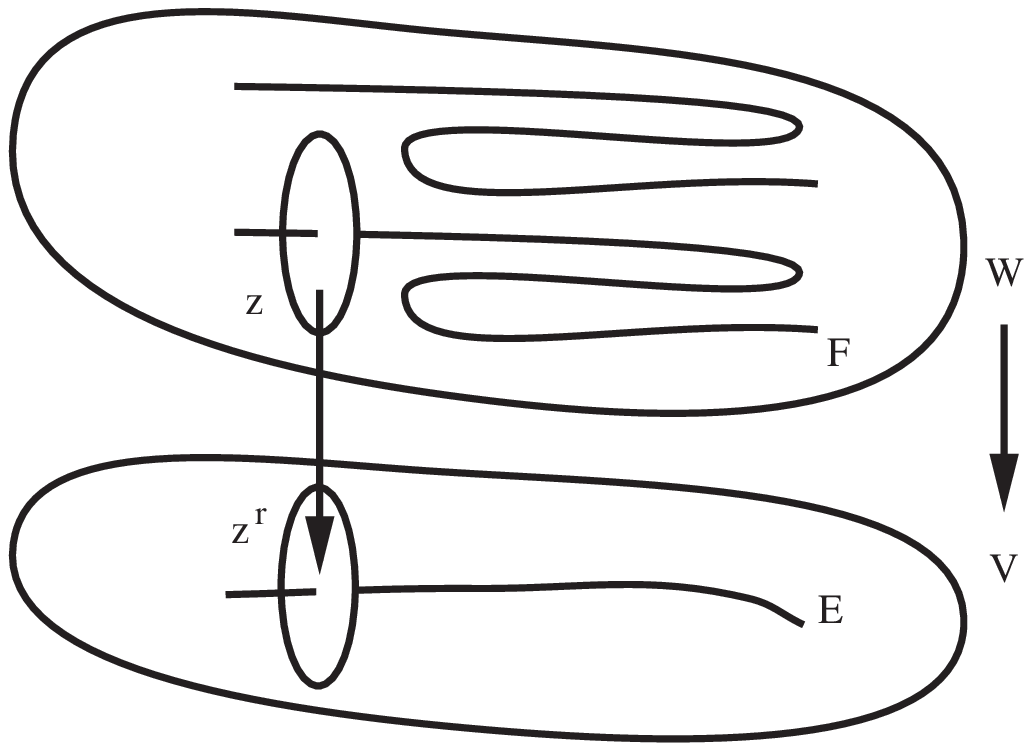}
\hfil
\par\vskip.2cm
\noindent
\hfil
{\bf Figure~2.} The cover $W\to V$ over the generic point of $E$.
\hfil
\endinsert

Thus we can assume that $E\subset V\to X$ is a prime divisor on a partial
resolution of $X$, and that $W$ is a blowup of $\C^n$ that dominates $V$:
 $$
\spreadmatrixlines{3pt}
\matrix
W&\longrightarrow\kern-0.6em&\C^n\\
\downarrow&&\downarrow\\
V&\to&X=\C^n/G
\endmatrix
 $$
(see Figure~2). Write $\sF=\{F\}$ for the set of all prime divisors
$F\subset W$ which dominate $E\subset V$. Under these assumptions, the
valuations extending $w=v_E$ are of the form $v_F$ where $F\in\sF$. Then $G$
acts transitively on $\sF$, and if $v=v_F$ then
 $$
\Stab v=\Stab F=\bigl\{g\in G\,\bigm|\,g(F)=F\bigr\}
 $$
and
 $$
\Ram v=\Ram F=\bigl\{g\in G\,\bigm|\,g\rest F=\id_F\bigr\}.
 $$
Note: the action of $G$ on $W$ is only birational; however, it is regular at
the generic point of every component $F$, and $g(F)=F$ and $g\rest F=\id_F$
are to be interpreted in this sense.

Now since $G$ acts transitively, the orbit space $G/\Stab F$ equals $\sF$.
Moreover, $\Stab F$ acts birationally on $F$ with birational quotient $E$,
so that $\Ram F\triangleleft\Stab F$ is a normal subgroup, and
 $$
\Stab F/\Ram F=\Gal(F/E)=\Gal(k(F)/k(E)).
 $$
Finally, $\Ram F$ is the cyclic group describing the behaviour of the
cover $W\to V$ along the general point of $F$. Since we are dealing there
with a codimension 1 submanifold of a complex manifold, the covering is
locally of the form $z\mapsto w=z^r$, so that $\Ram F=\Z/r$ is cyclic.

\subsubhead Ramification and discrepancy \endsubsubhead If $s$ is a canonical
form (that is, $n$@-form) on $V$ with $v_E(s)=m$ then
 $$
s=\text{(unit)}\cdot w^m\dd x_1\wedge\cdots\wedge\dd x_{n-1}\wedge\dd w
 $$
with $x_1,\dots,x_{n-1}$ local coordinates along $E$, and its pullback to
$W$ is of the form
 $$
\psi^*s=\text{(unit)}\cdot z^{mr+r-1}\dd x_1\wedge\cdots\wedge\dd
x_{n-1}\wedge\dd z,
 $$
so that $v_F(\psi^*s)=mr+r-1=rv_E(s)+r-1$. Therefore $a_E=\rth(a_F-(r-1))$.

As discussed in 2.3, a crucial point for us is that the discrepancy of $E$ is
already present in the cyclic quotient $\C^n/\Ram F$. In our picture, the
quotient $W/\Ram F$ contains a divisor $\overline F$, and the subsequent map
$W/\Ram F\to V$ is etale at a general point of $\overline F$.

 \subhead 2.7. Ramification for monomial valuations \endsubhead The definitions
of $\Stab v$ and $\Ram v$ for a discrete valuation $v$ of $k(\C^n)$ extend to
the whole of $\GL(n,\C)$ (without reference to a finite subgroup) as follows:
 $$
 \Stab'v=\Stab'F=\bigl\{g\in\GL(n,\C)\,\bigm|\,v\circ g=v\bigr\}
 =\bigl\{g\in\GL(n,\C)\,\bigm|\,g(F)=F\bigr\}
 $$
and
 $$
 \Ram'v
 =\bigl\{g\in\Stab'v\,\bigm|\,g^*=\id\text{ on }k(v)\bigr\}
 =\bigl\{g\in\Stab'F\,\bigm|\,g\rest F=\id_F\bigr\},
 $$
where $F\subset W$ is a prime divisor of a blowup $W\to\C^n$ such that
$v=v_F$. Here $\Stab'v=\Stab_{\GL}v$ denotes (temporarily) the stabiliser
of $v$ in $\GL$, and similarly for $\Ram'$.

We now determine these subgroups for monomial valuations. We first fix some
ideas and notation. Consider $\C^n$ with a fixed coordinate system
$x_1,\dots,x_n$. A primitive nonzero vector $\be=(b_1,\dots,b_n)\in\Z^n$ with
all $b_i\ge0$ corresponds to the $\C^\times$ action on $\C^n$ given by
 $$
\la\colon(x_1,\dots,x_n)\mapsto(\la^{b_1}x_1,\dots,\la^{b_n}x_n).
 $$
The action defines a decomposition of $\C^n$ as a direct sum of eigenspaces;
we simply group together the $x_i$ with the same $b_i$:
 $$
 \C^n=\bigoplus_{b\in\Z}S_b,
\quad\text{where}\quad S_b=\Spec k\bigl[x_i\,|\,b_i=b\bigr].\tag2.7.1
 $$
There is a corresponding filtration by increasing $b_i$:
 $$
T_0\subset\cdots\subset T_c\subset\cdots\subset\C^n,
\quad\text{where}\quad
T_c=\Spec k\bigl[x_i\,|\,b_i\le c\bigr]=\bigoplus_{b\le c}S_b,\tag2.7.2
 $$
and, of course, $\Gr_bT=T_b/T_{b-1}=S_b$.

 The weighting $\be=(b_1,\dots,b_n)$ defines a grading on
$k[x_1,\dots,x_n]$, and a weighted blowup
$\si_\be\colon W_\be\to\C^n$, whose exceptional locus $F_\be\subset W_\be$ is
$\Proj k[x_1,\dots,x_n]=\proj(b_1,\dots,b_n)$. Note that some of the $b_i$ may
be zero, and then $\Proj$ must be taken as a product of a conventional weighted
projective space with the affine space $S_0$, the centre of the blowup in
$\C^n$. The weighting $(b_1,\dots,b_n)$ also defines a monomial discrete
valuation $v_{F_\be}=v_\be$, taking $x_i\mapsto b_i$.

 \proclaim{Lemma} As above, $\be=(b_1,\dots,b_n)$ is a primitive nonzero vector
in the positive octant; suppose that $g\in\GL(n,\C)$.

 \roster
 \item $g\in\Stab'v_\be$ if and only if it preserves the filtration
$(2.7.2)$, that is, $g^*$ takes each $x_i$ to a linear combination of $x_j$
with\/ $b_j\le b_i$:
 $$
g^*(x_i)=\sum m_{ij}(g)x_j\quad\text{with $m_{ij}=0$ whenever $b_i>b_j$}.
\tag2.7.3
 $$

 \item Suppose that $g\in\Stab'v_\be$; then the action of $g$ on
$F_\be=\proj(b_1,\dots,b_n)$ is given by the induced map $\overline g$ on the
graded vector space $\bigoplus T_b/T_{b-1}$, that is,
 $$
\overline g^*(x_i)=\sum\overline m_{ij}(g)x_j\quad\text{with }
\overline m_{ij}=\cases
m_{ij}&\text{if\/ $b_j=b_i$,}\\
0&\text{if\/ $b_j<b_i$.}
\endcases
\tag2.7.4
 $$
In particular, $g\in\Ram'v_\be$ if and only if $\overline
g=\diag(\la^{b_1},\dots,\la^{b_n})$ for some $\la\in\C^\times$.

 \endroster
 \endproclaim

 \demo{Proof} This is really obvious, because the valuation $v_\be$ determines
the filtration (2.7.2). Indeed, the valuation of a linear form $\sum\al_ix_i$
for $v_\be$ is just the smallest $b_i$ for which $\al_i\ne0$. For (1),
$g\in\GL(n,\C)$ preserves the space of linear forms, and if it also preserves
$v_\be$ then it preserves the filtration (2.7.2).

If $g$ preserves the filtration (2.7.2), it also preserves the filtration
of $k[x_1,\dots,x_n]$ by the ideals of polynomials vanishing $k$ times along
$F_\be$:
 $$
I_k=\bigl\{f\,\bigm|\,v_\be(f)\ge k\bigr\}.
 $$
Therefore its action on $\C^n$ extends to a biregular action on the weighted
blowup $W_\be=\Proj\bigoplus_{k\ge0}I_k\to\C^n$. The restriction to
$F_\be=\proj(b_1,\dots,b_n)$ is in terms of the associated
$\Gr=\bigoplus_{k\ge0}I_k/I_{k+1}$, so the first part of (2) is clear. For
the final part, the condition $g\in\Ram'v_\be$ is that $g$ fixes $F_\be$
pointwise (on a dense open set). Thus this is the usual homogeneous
coordinate business: $(x_1:\cdots:x_n)=(y_1:\cdots:y_n)$ in
$\proj(b_1,\dots,b_n)$ if and only if $x_i=\la^{b_i}y_i$ for some
$\la\in\C^\times$. \QED\enddemo

 \proclaim{Corollary} Let $G\subset\GL(n,\C)$ be a finite group, and $g\in G$
an element of order $r$, written as usual $g=\rth(b_1,\dots,b_n)$ in
eigencoordinates $x_1,\dots,x_n$. Assume that $\be=(b_1,\dots,b_n)$ is a
primitive vector of $\Z^n$, and let $v=v_\be$ be the corresponding monomial
valuation.

Then an element $h\in G$ is in $\Stab v_\be$ if and only if it preserves the
direct sum decomposition $(2.7.1)$. Moreover, $h\in\Ram v_\be$ if and only
if $h=\diag(\ep^{b_1},\dots,\ep^{b_n})$, where $\ep$ is a root of $1$.
 \endproclaim

 \demo{Proof} Choose a $G$-invariant Hermitian metric on $\C^n$, so that all
the elements of $G$ are unitary. Then the direct sum (2.7.1) is the eigenspace
decomposition of the unitary operator $g$, and is thus orthogonal. Moreover,
$h$ is also a unitary operator, so it preserves the filtration (2.7.2) if and
only if it preserves the direct sum decomposition (2.7.1). The final part
follows directly from (2) of the lemma. \QED\enddemo

\subhead 2.8. Conclusion of proof of Theorem~1.4 \endsubhead We divide up the
proof into 5 easy steps.

\demo{Step~1. Construction of $g\mapsto v_g$} As we hinted above, the
correspondence
 $$
\Ga_1\to\bigl\{\text{valuations of $k(X)$}\bigr\},
 $$
sends $g$ to the monomial valuation $v_g$ of $k(X)$ corresponding to $g$. That
is, as described in 2.2, write $g=\rth(a_1,\dots,a_n)\in\Box$, to mean
that $g$ has order exactly $r$, and is
$\ep\mapsto\diag(\ep^{a_1},\dots,\ep^{a_n})$ in eigencoordinates
$x_1,\dots,x_n$. Now set $\wave v_g=v_\be$, where $\be=(a_1,\dots,a_n)$, and
 $$
v_g=\rth(\wave v_g{}\rest{k(X)})\colon k(X)\to\Z\cup\{\infty\}.
 $$
It's clear from the definition of discrete valuation that $v_g$ is (possibly a
multiple of) a discrete valuation of $k(X)$. \enddemo

\demo{Step~2. Ramification and discrepancy} We claim that for $g\in\Ga_1$, the
ramification group of $\wave v_g=v_\be$ in $G$ is exactly $\Span g$, and the
ramification degree is $r$. It follows from this that $v_g$ is in fact a
discrete valuation of $k(X)$, and, by the calculation in 2.4, that its
discrepancy is $\age g-1=0$.

Since $g=\rth(a_1,\dots,a_n)$, Corollary~2.7 says that the ramification
group of $\wave v_g$ is the subgroup of $G$ of elements of the form
$h=\diag(\ep^{a_1},\dots,\ep^{a_n})$, where $\ep$ is a root of 1. If the
order of $\ep$ divides $r$ then $h\in\Span g$, which is what we want. If not,
replacing $h$ by a suitable combination of $g,h$, we can assume that $\ep$
is a root of unity  of order larger than $r$; but then, because $\sum a_i=r$,
such an $h$ is not in $\SL(n,\C)$. \enddemo

\demo{Step~3. $v_g$ is well defined on conjugacy classes} If $g_1=hgh\1$, and
$x_1,\dots,x_n$ are eigencoordinates for $g$ then $h^*x_i$ are eigencoordinates
for $g_1$ with the same eigenvalues. Since $\wave v_g$ and $\wave v_{g_1}$ are
constructed as monomial valuations with respect to these eigencoordinates, it
follows at once that $\wave v_g\circ h=\wave v_{g_1}$, and therefore they
restrict to the same valuations of $k(X)=k(\C^n)^G$.

 In geometric terms, the weighted blowup of $\C^n$ corresponding to $g$
and $g_1$ fit together into a diagram of morphisms:
 $$
\spreadmatrixlines{3pt}
\matrix W_g&@>h>>&W_{g_1}\\
\downarrow&&\downarrow\\
\C^n&@>>h>&\C^n
\endmatrix
 $$
That is, $h$ takes the exceptional prime divisor $F_g$ into $F_{g_1}$, so that
they both map birationally to the same exceptional divisor of $X$. \enddemo

\demo{Step~4. Injectivity} This is the converse of Step~3: if $g,g'\in\Ga_1$
and $v_g=v_{g'}$ then $g$ and $g'$ are conjugate in $G$. The assumption is that
$\wave v_g$ and $\wave v_{g'}$ restrict to the same valuation of $k(X)$. Recall
that the Galois group $G=\Gal(k(\C^n)/X)$ acts transitively on the set of
extensions of a valuation of $k(X)$. Therefore $\wave v_{g'}=\wave v_g\circ h$
for some $h\in G$. We've seen in Step~3 that $\wave v_g\circ h=\wave v_{g_1}$
where $g_1=hgh\1$. Now $\wave v_{g_1}=\wave v_{g'}$ clearly implies that
$g'=g_1$. (Because they both induce the same filtration of the vector space
of linear forms, and give $x_i$ the same weights $b_i$.)
\enddemo

\demo{Step~5. Surjectivity: every crepant divisor of $X$ comes from this
construction} If $v_E$ is a crepant valuation of $k(X)$ and $\wave v=v_F$ a
valuation of $k(\C^n)$ extending $v$ then $\Ram\wave v$ is a cyclic group of
order $r$, the ramification degree. The discrepancy of $E$ depends only on $r$
and that of $F$, and is given by
 $$
a_E=\rth(a_F-(r-1)).
 $$
As explained in 2.3, only monomial valuations of $\C^n$ with nontrivial
ramification groups can give crepant valuations of $X$. These are all of the
form $v_g$ with $g\in\Ga$, and only those with $g\in\Ga_1$ are crepant.

This completes the proof of Theorem~1.4. \QED\enddemo

 \remark{Remark} One can also consider monomial valuations $v_g$ for conjugacy
classes $g\in\Ga_i/G$ with $i\ge2$. To get the right discrepancy
$a_{E_g}=\age-1$, and to keep the injectivity proved in Step~3, we restrict to
{\it primitive} elements $g$ (corresponding to a primitive vector
$g=\rth(a_1,\dots,a_n)\in\Box$). It's an exercise to generalise Steps~2--4
to these valuations. For example, the cyclic group $\frac14(1,1,1,1)$ (the
Gorenstein cone on the 4-fold Veronese embedding of $\proj^3$) has only one
crepant divisor. \endremark

\subhead 2.9. Proof of Corollary~1.5 \endsubhead This is clear. If a minimal
model $Y$ exists then the exceptional divisors occurring on it are the crepant
divisors and no others. Next, $H^2(Y,\Q)$ is spanned by algebraic cycles, the
first Chern classes of divisors (just use the exponential sequence and
$R^1f_*\Oh_Y=0$). But every divisor $D$ on $Y$ pushes down to divisor on $X$
some multiple of which is linearly equivalent to zero, because a quotient
singularity is analytically $\Q$@-factorial. Therefore some multiple of $D$ is
linearly equivalent to an exceptional divisor. \QED

\subhead 2.10. Proof of Theorem~1.6 \endsubhead For
$g=\rth(a_1,\dots,a_n)\in\Box$, the fixed locus $\Fix g\subset\C^n$ is the
linear subspace $S_0$ corresponding to the coordinates $x_i$ with $a_i=0$. In
particular, 0 is an isolated point of $\Fix g$ if and only if all $a_i>0$, or
equivalently, $g$ is in the interior of $\Box$. Thus we set
 $$
\Ga_1^{(0)}=\bigl\{g\,\bigm|\,\Fix g=\{0\}\bigr\}.
 $$

Note in passing that if $n=3$ and $g=\rth(a_1,a_2,a_3)$ then $a_1=0$ implies
$a_2+a_3=r$, so that
 $$
\Fix g\ne\{0\}\implies\text{$g$ is junior.}
 $$

We now restrict attention to the case $n=3$. Then $g\in\Ga_1^{(0)}$ is
equivalent to $g\1\in\Ga_2$: in eigencoordinates, if $g=\rth(a_1,a_2,a_3)$
with all $a_i>0$ then
 $$
g\1=(1,1,1)-\rth(a_1,a_2,a_3).
 $$
Thus $g\mapsto g\1$ gives the bijection $\Ga_1^{(0)}\to\Ga_2$, which is
clearly compatible with conjugacy. This proves the first bijection in
Theorem~1.6.

The second bijection is clear, since a divisor $E$ has compactly supported
class in $H^2$ if and only if $E$ contracts to $0\in X$. The third bijection is
just Poincar\'e duality. \QED

\heading 3. Examples \endheading

We start with the examples of the Klein quotient singularities or Du Val
surface
singularities, and explain why our result Theorem~1.4 gives the exceptional
components {\it together with their configuration}. It is perhaps remarkable
that our methods lead naturally to a ``McKay quiver'' structure on the
conjugacy classes of $G$ (rather than irreducible representations).

The Du Val singularity $D_{n+2}$ is the quotient $\C^2/\BD_{4n}$, where
$\BD_{4n}$ is the binary dihedral group of order $4n$ generated by
 $$
 A=\left(\matrix\ep&0\\0&\ep\1\endmatrix\right)
\quad\text{and}\quad
 B=\left(\matrix0&1\\-1&0\endmatrix\right),
 $$
and $\ep$ is a primitive $2n$th root of $1$. Thus $A^n=B^2=-1$, and
the elements of $\BD_{4n}$ are the $2n$ ``rotations'' $A^i$, and the
$2n$ ``reflections'' $A^iB=BA^{2n-i}$.

Conjugacy in $G$ gives $A^i\sim A^{2n-i}$, so there are $n$ conjugacy classes
of nonidentity rotations; and $ABA\1=A^2B$, so that $B\sim A^{2i}B$ and $AB\sim
A^{2i+1}B$. The effect of conjugacy on the cyclic subgroups
$\Span{A^iB}\iso\Z/4$ depends on $n$ even or odd: if $n$ is even then $\Span B$
and $\Span{AB}$ are not conjugate, but $B\sim A^nB=-B=B^3$ and
$AB\sim(AB)^3=BA$. If $n$ is odd then all the subgroups $\Span{A^iB}$ are
conjugate, but the elements of order 4 within any subgroup are not conjugate.

\define\myc#1{\kern-1cm#1\kern-1cm}    % to centre a matrix entry
\define\hedge{\kern1mm{\qquad\over}\kern1mm} % horizontal edge
\define\vedge{\vrule height20pt depth0pt width0.4pt
 \vrule height 22pt depth0pt width0pt} % vertical edge

 \example{Example 1 \rm (Binary dihedral group $\BD_8$, singularity $D_4$)} In
this case the two matrixes are
 $$
A=\pmatrix i&0\\0&-i\endpmatrix\quad\text{and}\quad
B=\pmatrix0&1\\-1&0\endpmatrix \in\SL(2,\C),
 $$
where $i=\sqrt{-1}$. Thus $G$ has three cyclic subgroups of order 4,
$\Span A$, $\Span B$ and $\Span{AB}$, with $A^2=B^2=(AB)^2=-1$, which we
draw as follows:

\smallskip
 $$
\matrix
\format \r\qquad&&\c\\
&\myc A&&\myc{A^2}&&\myc{A^3}\\
\Span A:&\bigcirc&\hedge&\bigcirc&\hedge&\bigcirc\\
\vspace{.5cm}&\myc B&&\myc{B^2}&&\myc{B^3}\\
\Span B:&\bigcirc&\hedge&\bigcirc&\hedge&\bigcirc\\
\vspace{.5cm}&\myc{AB}&&\myc{(AB)^2}&&\myc{(AB)^3}\\
\Span{AB}:&\bigcirc&\hedge&\bigcirc&\hedge&\bigcirc\\ \endmatrix
\tag1
 $$

\smallskip The pictures are the exceptional divisors in the minimal resolution
of $\C^2/\Span g$ as given by toric geometry. The nodes represent elements of
$\Span g$, and they are joined by an edge if they are neighbouring points in
the
junior simplex (interval) $\De_1$, that is, if the corresponding curves
intersect on the minimal resolution.

Conjugacy in $G$ gives
 $$
A\sim A^3,\quad B\sim B^3,\quad AB\sim(AB)^3.
 $$
Thus identifying equal and conjugate group elements gives the familiar $D_4$
diagram:

 $$
\matrix
\format &\c\\
\myc A&&\myc{-1}&&\myc B\\
\bigcirc&\hedge&\bigcirc&\hedge&\bigcirc\\
&&\vedge\\
&\myc{AB}&\bigcirc
\endmatrix\tag2
 $$

\endexample

\example{Example 2 \rm (Binary dihedral group $\BD_{8l+4}$, singularity
$D_{2l+3}$)} In this case $A$ has order $4l+2$. We should draw all the maximal
cyclic subgroups, then consider the effect of conjugacy. As noted above, the
cyclic subgroups $\Span{A^iB}$ are all conjugate, so we only draw $\Span A$
and $\Span B$:

\smallskip
 $$
 \matrix
\format \r\qquad&&\c\\
&\myc A&&\myc{A^2}&&\myc{A^{4l}}&&\myc{A^{4l+1}}\\
\Span A:&\bigcirc&\hedge&\bigcirc&\quad\cdots\quad&\bigcirc&\hedge&\bigcirc\\
\vspace{.5cm}&\myc B&&\myc{B^2}&&\myc{B^3}\\
\Span B:&\bigcirc&\hedge&\bigcirc&\hedge&\bigcirc\\
 \endmatrix
\tag3
 $$

\smallskip Now $A^{2l+1}=B^2=-1$, and conjugacy gives $A^i\sim A^{4l+2-i}$, so
that identifying equal and conjugate elements gives
 $$
\matrix
\format &\c\\
\myc B&&\myc{-1}&&\myc{A^{2l}}&&\myc{A^2}&&\myc A\\
\bigcirc&\hedge&\bigcirc&\hedge&\bigcirc&\quad\cdots\quad&
\bigcirc&\hedge&\bigcirc\\
&&\vedge\\
&\myc{B^3}&\bigcirc
\endmatrix\tag4
 $$

The case of the binary dihedral group $\BD_{8l}$, giving the singularity
$D_{2l+2}$ is similar. The only difference is that to see the action of
conjugacy, we have to draw both $\Span B$ and $\Span{AB}$, and then fold each
of them in two.

\endexample

\example{Example 3 \rm (Binary tetrahedral group $\BT_{48}$, singularity
$E_6$)}
The singularity $E_6$ is the quotient of $\C^2$ by the group $\BT_{48}$
generated by the three matrixes
 $$
 A=\pmatrix i&0\\0&-i\endpmatrix,\quad
 B=\pmatrix0&1\\-1&0\endpmatrix,\quad\text{and}\quad
 C=\frac1{\sqrt2}\pmatrix e^7&e^7\\e^5&e\endpmatrix\in\SL(2,\C),
 $$
where $e$ is a primitive 8th root of 1 (for example, $e={1+i\over\sqrt2}$).
Thus
 $$
\BT_{48}=\bigl<A,B,C\bigm|
A^2=B^2=C^3=-1,\;BA^iB^{-1}=A^{-i},\;CAC^{-1}=B^3\bigr>.
 $$
We draw the cyclic subgroups $\Span A$ and $\Span C$ of $\BT_{48}$ following
the
previous recipe:

\smallskip
 $$
\matrix
\format \r\qquad&&\c\\
&\myc A&&\myc{A^2}&&\myc{A^3}\\
\Span A:&\bigcirc&\hedge&\bigcirc&\hedge&\bigcirc\\
\vspace{.5cm}&\myc C&&\myc{C^2}&&\myc{C^3}&&\myc{C^4}&&\myc{C^5}\\
\Span C:&\bigcirc&\hedge&\bigcirc&\hedge
&\bigcirc&\hedge&\bigcirc&\hedge&\bigcirc\\
\endmatrix
\tag5
 $$
\smallskip

It is easy to see the effect of conjugacy in $\BT_{48}$:
 $$
A\sim A^3\sim B\sim B^3\sim AB\sim BA,\quad\text{and}\quad A^kC^j\sim C^j
\text{ for all $j,k$},
 $$
and the elements of $\Span C$ are not conjugate to one another. Thus
identifying
equal and conjugate elements of these two groups gives the resolution graph of
$E_6$:
 $$
\matrix
\format &\c\\
\myc C&&\myc{C^2}&&\myc{-1}&&\myc{C^4}&&\myc{C^5}\\
\bigcirc&\hedge&\bigcirc&\hedge&\bigcirc&\hedge&\bigcirc&\hedge&\bigcirc\\
&&&&\vedge\\
&&&\myc{A}&\bigcirc
\endmatrix
\tag6
 $$

\endexample

Now we see some 3-dimensional examples. (More examples are furnished by the
crepant resolutions of $\C^3/G$ of \cite{Roan}, \cite{Ito1--3} and
\cite{Markushevich}.)

\example{Example 4 \rm (Trihedral group)} A {\it trihedral} group is a subgroup
$G$ generated by a diagonal subgroup $H\subset\SL(3,\C)$, together with the
cyclic permutation matrix
 $$
T=\pmatrix0&1&0\\0&0&1\\1&0&0\endpmatrix.
 $$
The Abelian subgroup $H$ is generated by one or more diagonal matrixes
 $$
A_i=\frac1{r_i}(a_i,b_i,c_i),\quad\text{with}\quad a_i+b_i+c_i=r_i,
 $$
meaning as usual that $A_i=\diag(\ep^{a_i},\ep^{b_i},\ep^{c_i})$ where $\ep$ is
a primitive $r_i$th root of 1. We can of course assume that $H$ has cyclic
symmetry, and then it is a normal subgroup $H\triangleleft G$. The quotient of
$\C^3$ by $G$ is a {\it trihedral quotient singularity}. A crepant resolution
of singularities $Y\to X$ is described in \cite{Ito1--2}.

We work out the special case
$H=\Span{\third(0,1,2),\third(1,2,0),\third(2,0,1)}$. Obviously,
$G=\Span{\third(0,1,2),T}$, and is a group of order 27.

It turns out that $G$ has only one senior element $\third(2,2,2)$, and 25
juniors. Note that, as we saw in 2.10, for subgroups of $\SL(3,\C)$, all
elements with nonisolated fixed points are junior; this applies to every
element of $G$ except $\third(1,1,1)$ and $\third(2,2,2)$ because, for
example, $\Fix\third(0,1,2)=x$@-axis, and every element $hT$ is obviously
conjugate in $\GL(3,\C)$ to $T$, which fixes the line
$\{(x,x,x)\}\subset\C^3$. Thus

$$
\spreadmatrixlines{3pt}
\Ga_1=
\left\{
\matrix
\third(1,1,1),\;
\third(0,1,2),\;
\third(1,2,0),\;
\third(2,0,1),\\
\third(0,2,1),\;
\third(2,1,0),\;
\third(1,0,2)
\endmatrix
\right\} \cup HT \cup HT^2
$$

Under conjugacy in $G$,
 $$
\third(a,b,c)\sim\third(b,c,a),\quad
\text{for all $a,b,c\in\{0,1,2\}$,}
 $$
and
 $$
T\sim\third(a,b,c)T,\quad\text{and}\quad T^2\sim\third(a,b,c)T^2
\quad\text{for all $a\ne b$.}
 $$
Thus the exceptional divisors in the minimal resolution $Y\to X=\C^3/G$
correspond one-to-one with the following junior conjugacy classes:
 $$
\spreadmatrixlines{4pt}
\matrix
\third(0,1,2),\;\third(0,2,1),\;\third(1,1,1),\\
T,\;\third(1,1,1)T,\;\third(2,2,2)T,\\
T^2,\;\third(1,1,1)T^2,\;\third(2,2,2)T^2.
\endmatrix
 $$

These elements correspond to a basis of $H^2(Y,\Q)$, so that $h^2(Y,\Q)=9$. On
the other hand, by Corollary~1.6, $H^4(Y,\Q)$ has a basis corresponding to
the single element $\Ga_2=\{\third(2,2,2)\}$, and its Poincar\'e dual
$H_c^2(Y,\Q)$ has a basis corresponding to $\Ga_1^{(0)}=\{\third(1,1,1)\}$.
In particular the Euler number $e(Y)$ is given by
 $$
\align
e(Y)&{}=h^0(Y,\Q)+h^2(Y,\Q)+h^4(Y,\Q)=1+9+1=11\\
 &{}=\sharp\{\text{conjugacy class in $G$}\}.
\endalign
 $$

\endexample

\example{Example 5 \rm (Icosahedral group)} $G_{60}\subset\SO(3)$ is the
group of rotations of the icosahedron, acting on $\C^3$ by the inclusion
$\SO(3)\subset\SL(3)$. Generating matrixes and defining relations are given
for example in \cite{Miller--Blichfeldt--Dickson}, p.~250. The reader
may eke out our poor words of explanation by imagining the picture of the
icosahedron.

It is easy to see that $G_{60}$ has just 5 conjugacy classes:

\roster
\item the identity;

\item a rotation through $\pi$ about an axis through the midpoint of an edge;

\item a rotation through $2\pi/3$ about an axis through the midpoint of a face;

\item a rotation through $2\pi/5$ about an axis through a vertex;

\item a rotation through $4\pi/5$ about the same axis.
\endroster

Each rotation fixes an axis, and is therefore junior by the remark in 2.10.
The quotient $\C^3/G_{60}$ has 3 curves of transversal Du Val singularities:
a curve of $A_1$ under the axis of (2); a curve of $A_2$ under the axis of
(3), and a curve $C$ of transversal $A_4$ singularities under the common
axis of (4) and (5). Resolving these curves gives 4 crepant divisors. For
example, blowing up the curve $C$ of $A_4$ singularities
gives two surfaces which (outside the origin) are $\proj^1$ bundles over
a curve that is a double cover of $C$.

\endexample

\heading 4. Discussion \endheading

\subheading{4.1} Theorem~1.6 implies directly Conjecture~1.1 for $n=3$, with
the proviso that you have to fix roots of 1 to get a canonical correspondence.
To see that $Y$ cannot have any odd-dimensional homology, it is enough to note
that there exists a resolution of $X$ dominated by a variety $W$ obtained by
a sequence of blow-ups of $\C^n$ in linear subspaces. Then $H^*(W)$ has no
odd-dimensional cohomology, and $H^*(Y)\subset H^*(W)$.

\subheading{4.2} The minimal model is in fact nonsingular: indeed, by
standard results of Mori theory, although the minimal model is usually not
unique, all minimal models have the same local singularities. Moreover,
Roan and others have proved that a smooth crepant resolution exists in all
cases. Unfortunately, the present abstract point of view does not (as yet)
offer any clues as to why this should be so.

\subheading{4.3} There should be a theoretical proof that the minimal model
$Y\to X$ is nonsingular, even though we are not clever enough to find it
at present. For example, if $Y$ has factorial cDV points, then these have
local deformations. It should follow from Friedman's local-to-global spectral
sequence that $Y$ has deformations inducing all possible deformations of the
local singularities, and in particular, should have a smoothing $Y_t$; then
$Y_t$ also contracts to a variety $X_t$. Alternatively the Milnor fibre shows
that $Y_t$ must have a middle cohomology of Hodge type $H^{1,2}\oplus
H^{2,1}$---but so what? Or maybe the singularities would have homotopy
properties that contradict that the link of $X$ is 2-connected.

\subheading{4.4} The shape of Theorems~1.4 and 1.6 is very suggestive: the
canonical basis for $H^2(Y,\Z)$ involves junior conjugacy classes of
$\Ga=G(-1)$, which is probably the only thing that makes sense in terms of
motivic weights. The basis of $H^4(Y,\Q)$ is dual to a basis of
$H^2_c(Y,\Q)$, so to take account of the roots of 1 correctly, you
presumably have to take $\Hom$ into $\Z(-3)$, and end up with another Tate
twist of the set of conjugacy classes $\Ga_2/G$ before getting a canonical
identification with the basis of $H^4(Y,\Q)$ (?).

\subheading{4.5. Crepant resolutions in higher dimension} In dimension
$\ge4$, when does $X=\C^n/G$ have a crepant resolution $Y\to X$? The known
counterexamples are the terminal quotient singularities, for example,
$\rth(i,r-i,j,r-j)$ (see also \cite{Mori--Morrison--Morrison}). For these,
there are not enough junior elements: the unit cube $\Box$ in $\R^4$
contains $r-1$ elements on the $\Ga_2$ hyperplane $\sum\al_i=2$, and no
junior elements at all. Thus there are no crepant divisors to pull out.

It seems plausible that the following 2 conditions to the effect that $G$ has
lots of junior elements are necessary conditions: for each {\it maximal}
cyclic subgroup $\Z/r\subset G$, choose eigencoordinates and write
$L=\Z^n+\Z\cdot\rth(a_1,\dots,a_n)$ as usual. Then

(i) every point of $L\cap\Box$ is a positive integral combination of
points of $L\cap\De_1$.

(ii) the positive octant of $L$ has a subdivision into basic cones with
1-skeleton contained in $\De_1$.

Obviously (ii) implies (i), and, for Abelian groups, by toric geometry,
(ii) is a necessary and sufficient condition.

\subheading{4.6. Multiplicative structure} Assume that a crepant resolution
$Y\to X$ exists, and that $n=4$ or $5$, to fix ideas. There is one respect
in which Conjecture~1.1 looks bad: to get a canonical correspondence, as we
have seen, you have to introduce an age grading, and this will only be
canonical if you twist $G$ to $G(-1)$. No other twist has a canonical age
grading. But then the theory of weights in the philosophy of motives make it
ungrammatical to put elements of $G(-1)$ in correspondence with elements of
$H^4$.

One way of trying to rescue the conjecture might be to put elements of
$H^4$ in correspondence with quadratic monomials in $\Ga_1$, maybe with
relations between them deduced from a multiplicative structure; this seems
to be good for Abelian quotients (compare \cite{Batyrev--Dais}). Some natural
questions arise in this connection. Presumably it is true that the
cohomology ring $H^*(Y,\Q)$ is generated by $H^2$? Is there a reasonable
multiplicative structure involving the conjugacy classes of $G$?

\subheading{4.7. Remark} The easy arguments in this paper have left to
one side the most important aspect of the McKay correspondence for
$\SL(2,\C)$, namely, the relation with irreducible representations.
Understanding this point might begin to explain a lot of things. Some
indications:

 (1) In 2.5 we needed to find an extraction $V\to X$ pulling out exactly
the divisor $E$ corresponding to a conjugacy class, by analogy with the
weighted blowup of 2.3 in the cyclic case. This could be provided by
quiver-theoretic methods: view $X=\C^n/G$ as a moduli space of
representations of a quiver (e.g.\ the McKay quiver of $G$). Then $X$ is a
categorical quotient of a big vector space by a big algebraic group, and
birational modifications of $X$ are naturally obtained by ``varying the
linearisation''.

 (2) The age of $g\in\Ga$ is defined as the elementary symmetric function $\age
g=\rth\sum a_i$. Do the other elementary symmetric functions (e.g.\
$\frac1{r^2}\sum a_ia_j$) also have some significance for the McKay
correspondence?

\heading References \endheading

\parindent=0pt

V. Batyrev, Strong McKay correspondence, lecture at Trento conference, Jun 1994

Victor V. Batyrev and Dimitrios I. Dais, Strong McKay correspondence,
String-theoretic Hodge Numbers and Mirror Symmetry, Duke server
alg-geom/9410001, to appear in Topology {\bf35} (1996)

V. Danilov, The geometry of toric varieties, Uspekhi Mat.\ Nauk {\bf
33}:2 (1978) 85--134 = Russian Math.\ Surveys {\bf 33}:2 (1978) 97--154

V. Danilov, The birational geometry of toric $3$@-folds, Izv.\
Akad.\ Nauk SSSR Ser.\ Mat.\ {\bf 46} (1982), 972--981 = Math.\ USSR -- Izv.\
{\bf 21} (1983)

L. Dixon, J. Harvey, C. Vafa and E. Witten, Strings on orbifolds, Nucl.\ Phys.\
B {\bf261} (1985) 678--686. II, Nucl.\ Phys.\ B {\bf274} (1986) 285--314

G. Gonz\'alez-Sprinberg and J.-L. Verdier, Points doubles rationnels et
repr\'esent\-ations de groupes, C.~R. Acad. Sci. Paris, S\'er.~I, Math.\
{\bf 293} (1981), 111--113

G. Gonz\'alez-Sprinberg and J.-L. Verdier, Construction g\'eometrique
de la correspondance de McKay, Ann. Sci. \'Ecole Norm. Sup. (4) {\bf16}
(1983), 409--449

F. Hirzebruch and T. H\"ofer, On the Euler number of an orbifold, Math.\ Ann.\
{\bf286} (1990), 255--260

Y. Ito, Crepant resolutions of trihedral singularities, Proc.\ Japan Acad.\
{\bf70} (1994), 131--136

Y. Ito, Crepant resolution of trihedral singularities and the orbifold Euler
characteristic, Internat.\ J. Math.\ {\bf6}:1, 1995, 33--43

Y. Ito, Gorenstein quotient singularities of monomial type in dimension three,
J.\ Math.\ Sci.\ Univ.\ Tokyo {\bf2} (1995), 419--440. Duke server
alg-geom/9406001

J. Koll\'ar, Flips, flops etc., in Surveys in differential geometry, {\bf1}
(1991) 113--199

D. Markushevich, Resolution of $\C^3/H_{168}$, Technion preprint, 1993

G. A. Miller, H. F. Blichfeldt and L. E. Dickson, Theory and
applications of finite groups, Wiley, London, 1917, reprinted Dover 1961

J. McKay, Graphs, singularities, and finite group, in Santa Cruz
conference on finite groups (Santa Cruz, 1979), Proc.\ Symp.\ Pure Math.\ {\bf
37}, AMS (1980), 183--186

S. Mori, D. R. Morrison and I. Morrison, On four-dimensional terminal
singularities, Math.\ Comp.\ {\bf51} (1988), 769--786

S-S. Roan, Minimal resolutions of Gorenstein orbifolds in dimension three,
Preprint R940606--1, Acad.\ Sinica, Taibei, Jun 1994, to appear in Topology
{\bf35} (1996)

[C3-f] M. Reid, Canonical 3-folds, in G\'eom\'etrie alg\'ebrique, (Angers 1979)
(A. Beauville, Ed.), Sijthoff and Noordhoof, Alphen aan den Rijn, 1980,
273--310

[YPG] M. Reid, Young person's guide to canonical singularities, in Algebraic
geometry (Bowdoin 1985), Proc.\ Symp.\ Pure Math.\ {\bf 46}, AMS, 1987,
345--416

[Z--S] O. Zariski and P. Samuel, Commutative Algebra, Vols.~I and~II, Van
Nostrand, 1960

\enddocument
\bye